\documentclass[11pt]{article}
\usepackage{latexsym}
\usepackage{amssymb}
\usepackage{amsmath,amscd}
\usepackage{amsthm}
\usepackage{mathrsfs}
\allowdisplaybreaks

\newtheorem{theorem}{Theorem}[section]
\newtheorem{lemma}{Lemma}[section]
\newtheorem{proposition}{Proposition}[section]
\newtheorem{definition}{Definition}[section]
\newtheorem{remark}{Remark}[section]
\newtheorem{corollary}{Corollary}[section]
\newtheorem{example}{Example}[section]

\begin{document}

\title{\bf The random case of Conley's
theorem: II. The complete Lyapunov function\thanks{Published in:
{\em Nonlinearity} {\bf 20} (2007), 1017--1030.}}

\author{Zhenxin Liu\\
{\small College of Mathematics, and Key Laboratory of Symbolic
Computation and Knowledge} \\ {\small  Engineering of Ministry of
Education, Jilin University,  Changchun 130012,
P.R. China}\\
{\small E-mail: zxliu@jlu.edu.cn}}

\date{}
\maketitle

\begin{abstract}
 Conley in \cite{Con} constructed a complete Lyapunov function for a
 flow on compact metric space which is constant on orbits in the chain recurrent set and is
strictly decreasing on orbits outside the chain recurrent set. This
indicates that the dynamical complexity focuses on the chain
recurrent set and the dynamical behavior outside the chain recurrent
set is quite simple. In this paper, a similar result is obtained for
random dynamical systems under the assumption that the base space
$(\Omega,\mathcal F,\mathbb P)$ is a separable metric space endowed
with a probability measure. By constructing a complete Lyapunov
function, which is constant on orbits in the random chain recurrent
set and is strictly decreasing on orbits outside the random chain
recurrent set, the random case of Conley's fundamental theorem of
dynamical systems is obtained. Furthermore, this result for random
dynamical systems is generalized to noncompact state spaces. 
\end{abstract}

\maketitle

\section{Introduction and main result}

In \cite{Con}, Conley presented an important result in the study of
dynamical systems--Conley decomposition theorem. This result, named
after him, is also called ``fundamental theorem of dynamical
systems" \cite{Nor}, which can be stated very simply:

\begin{theorem}\label{Conley}{\bf (Conley's fundamental theorem of dynamical systems).}
Any flow on a compact metric space decomposes the space into a chain
recurrent part and a gradient-like part.
\end{theorem}

The proof of the theorem relies on two important steps: the
characterization of chain recurrent set in terms of attractors and
the construction of complete Lyapunov function which provides the
gradient-like structure for the theorem. To be more specific, assume
$(X,d_X)$ is a compact metric space and $\varphi$ is a flow on $X$.
For given $\epsilon, T>0$, a finite sequence $(x_0,t_0)$,
$(x_1,t_1)$, $\cdots$, $(x_n,t_n)$ in $X\times (0,\infty)$ is called
an $\epsilon$-$T$-chain for $\varphi$ if
\[
d_X(\varphi(t_j,x_j),x_{j+1})<\epsilon,~t_j\ge T
\]
for $j=0,1,\cdots,n-1$. Here $n$ is called the {\em length} of the
chain. A point $p\in X$ is called {\em chain recurrent} if for any
$\epsilon, T>0$, there is an $\epsilon$-$T$-chain with the length at
least $1$ which begins and ends at $p$. The set of all chain
recurrent points, denoted by $\mathcal{CR}(\varphi)$, is called the
{\em chain recurrent set} of $\varphi$. Conley \cite{Con} showed
that the chain recurrent set $\mathcal{CR}(\varphi)$ equals the
intersection of the union of attractor-repeller pair $(A,R)$ as $A$
varies over the collection of attractors of $\varphi$, i.e.
\[
\mathcal{CR}(\varphi)=\bigcap [A\cup R].
\]
Based on this characterization, he constructed a complete Lyapunov
function. A {\em complete Lyapunov function} for $\varphi$ on $X$ is
a continuous, real-valued function $L$ with properties: (1) it is
strictly decreasing on orbits outside the chain recurrent set while
is constant on orbits in the chain recurrent set; (2)
$L(\mathcal{CR}(\varphi))$ is nowhere dense; and (3) it separates
different components of the chain recurrent set. A dynamical system
is called {\em gradient-like} if there exists some continuous
real-valued function which is strictly decreasing on nonconstant
solutions, see \cite{Con} for details of these definitions and some
interesting examples.

Conley's result was adapted for maps on compact spaces by Franks
\cite{Fra}, was later established for maps on locally compact metric
spaces by Hurley \cite{Hu0,Hu1}, and was extended by Hurley
\cite{Hu2} for semiflows and maps on arbitrary metric spaces. To
extend Conley decomposition theorem to random dynamical systems
(RDS), Liu \cite{Liu} introduced the definitions of random chain
recurrent set and appropriate random attractor. He then took the
first step towards the random case of Conley decomposition
theorem--the characterization of random chain recurrent set in terms
of random attractors. In this paper, we intend to complete the
second step--constructing a complete Lyapunov function, which
provides the gradient-like structure for RDS. Here, ``random
gradient-like" could be understood in this way: the skew-product
flow $\Theta_t(\omega,x):=(\theta_t\omega,\varphi(t,\omega)x)$
corresponding to RDS $\varphi$ (see \cite{Ar1,Ar2,Chu} for details)
flows downhill with respect to a function (i.e. complete Lyapunov
function) along non-stationary solutions (a solution $u(\omega)$ is
called {\em stationary} if
$\varphi(t,\omega)u(\omega)=u(\theta_t\omega)$). The main result of
this paper states as follows:

\begin{theorem}\label{main}
Assume $(\Omega,\mathcal F,\mathbb P)$ is a separable metric space
endowed with a probability measure, where $\mathcal F=\mathcal
B(\Omega)$, the Borel $\sigma$-algebra on $\Omega$. Let $X$ be a
compact metric space, and $\varphi$ be an RDS with the base space
$(\Omega,\mathcal F,\mathbb P)$ and the state space $X$, then there
exists a complete Lyapunov function for $\varphi$ (defined in
Definition \ref{LY}) on $X$.
\end{theorem}

By the main theorem in \cite{Liu} and Theorem \ref{main} above we
arrive at the random case of Conley decomposition theorem:

\begin{theorem}{\bf (random Conley decomposition theorem).}
Any random dynamical system (with a separable metric space endowed
with a probability measure as base space) on a compact metric space
decomposes the space into a random chain recurrent part and a random
gradient-like part.
\end{theorem}

\section{Preliminaries}

In this section two preliminary definitions are given. For other
definitions, notations and preliminary propositions we refer the
readers to \cite{Liu} for details.

\begin{definition}\label{set}
Let $X$ be a metric space with a metric $d_X$. A set-valued map
$\omega\mapsto D(\omega)$ taking values in the closed/compact
subsets of $X$ is said to be a {\em random closed/compact set} if
the mapping $\omega\mapsto {\rm dist}_X(x,D(\omega))$ is measurable
for any $x \in X$, where ${\rm dist}_X(x,B):=\inf_{y\in B}d_X(x,y)$.
A set-valued map $\omega\mapsto U(\omega)$ taking values in the open
subsets of $X$ is said to be a {\em random open set} if
$\omega\mapsto U^c(\omega)$ is a random closed set, where $U^c$
denotes the complement of $U$.
\end{definition}

\begin{definition}
For any given random set $D(\omega)$, we denote $\Omega_D(\omega)$
the {\em omega-limit set of $D(\omega)$}, which is determined as
follows:
\[
\Omega_D(\omega):=\bigcap_{T\ge 0}\overline{\bigcup_{t\ge
T}\phi(t,\theta_{-t}\omega)D(\theta_{-t}\omega)};
\]
and we denote $\alpha_D(\omega)$ the {\em alpha-limit set of
$D(\omega)$}, which is determined as follows:
\[
\alpha_D(\omega):=\bigcap_{T\ge 0}\overline{\bigcup_{t\ge
T}\phi(-t,\theta_{t}\omega)D(\theta_{t}\omega)}.
\]
\end{definition}

In Sections 2--4, we will always assume that $X$ is a compact metric
space, therefore it is a {\em Polish space}, i.e. a separable
complete metric space.

\begin{remark}\label{att}\rm
Assume $U(\omega)$ is a random pre-attractor and $A(\omega)$ is the
random local attractor inside $U(\omega)$ (see \cite{Liu} for the
definitions), then it is easy to see that
\begin{align*}
A(\omega)&=\bigcap_{n\in\mathbb N}\overline {\bigcup_{s\ge
nT(\omega)}\varphi(s,\theta_{-s}\omega)U(\theta_{-s}\omega)}\\
&=\bigcap_{t\ge T(\omega)}\overline {\bigcup_{s\ge
t}\varphi(s,\theta_{-s}\omega)U(\theta_{-s}\omega)}\\
&=\bigcap_{t\ge 0}\overline {\bigcup_{s\ge
t}\varphi(s,\theta_{-s}\omega)U(\theta_{-s}\omega)}\\
&=\Omega_U(\omega).
\end{align*}
Conversely, if a random open set $U(\omega)$ satisfies
$\Omega_U(\omega)\subset U(\omega)$, we define
\[
T(\omega):=\inf\{t\in\mathbb
R^+|~\varphi(t,\theta_{-t}\omega)\overline{U(\theta_{-t}\omega)}\subset
V(\omega)\},
\]
where $V(\omega)\subset U(\omega)$ is a closed neighborhood of
$\Omega_U(\omega)$. Then, by a similar argument to the proof of
Lemma 3.5 in \cite{Liu}, we obtain that $T(\omega)$ is measurable.
Hence $U(\omega)$ is a random pre-attractor and $\Omega_U(\omega)$
is the random attractor determined by $U(\omega)$. That is, an
invariant random compact set is an attractor if and only if it is
the omega-limit set of one of its random neighborhoods.
\end{remark}

\section{Improvement of main result in \cite{Liu}}

\begin{lemma}\label{lem}
Assume $U(\omega)$ is a random open set and an invariant random
compact set $A(\omega)\subset U(\omega)$ satisfies that
$\Omega_U(\omega)=A(\omega)$, then there exists a forward invariant
random open set $\tilde U(\omega)$ with the same properties as
$U(\omega)$.
\end{lemma}
\noindent{\bf Proof.} Let
\[
\tilde U(\omega):={\rm int}[\overline{\bigcup_{t\ge
0}\varphi(t,\theta_{-t}\omega)U(\theta_{-t}\omega)}],
\]
then by Proposition 1.5.1 of \cite{Chu}, $\tilde U(\omega)$ is a
forward invariant random open set and $A(\omega)\subset\tilde
U(\omega)$ (note that $U(\omega)\subset\tilde U(\omega)$). Now we
show $\Omega_{\tilde U}(\omega)=A(\omega)$.
\begin{align}
\Omega_{\tilde U}(\omega)&=\bigcap_{T\ge 0}\overline{\bigcup_{s\ge
T}\varphi(s,\theta_{-s}\omega)\tilde U(\theta_{-s}\omega)}\nonumber\\
&=\bigcap_{T\ge 0}\overline{\bigcup_{s\ge
T}[\varphi(s,\theta_{-s}\omega)\bigcup_{t\ge
0}\varphi(t,\theta_{-t}\circ\theta_{-s}\omega)U(\theta_{-t}\circ\theta_{-s}\omega)]}\label{equa}\\
&=\bigcap_{T\ge 0}\overline{\bigcup_{s\ge T}\bigcup_{t\ge
0}\varphi(s,\theta_{-s}\omega)\circ\varphi(t,\theta_{-t}\circ\theta_{-s}\omega)U(\theta_{-t}\circ\theta_{-s}\omega)}
\nonumber\\
&=\bigcap_{T\ge 0}\overline{\bigcup_{s\ge T}\bigcup_{t\ge
0}\varphi(s+t,\theta_{-s-t}\omega)U(\theta_{-s-t}\omega)}
\nonumber\\
&=\bigcap_{T\ge 0}\overline{\bigcup_{s\ge
T}\varphi(s,\theta_{-s}\omega)U(\theta_{-s}\omega)} \nonumber\\
&=\Omega_U(\omega)=A(\omega), \nonumber
\end{align}
where (\ref{equa}) holds because for any random set $D(\omega)$ we
have
\[ \overline{\bigcup_{t\ge
T}\varphi(t,\theta_{-t}\omega)D(\theta_{-t}\omega)}=
\overline{\bigcup_{t\ge
T}\varphi(t,\theta_{-t}\omega)\overline{D(\theta_{-t}\omega)}}.
\]
This completes the proof. \hfill$\square$

In contrast to Lemma 3.3 in \cite{Liu}, we have the following finer
result.
\begin{corollary}\label{Ba}
For arbitrary random local attractor $A(\omega)$, the random basin
$B(A)(\omega)$ determined by $A(\omega)$ is an invariant random open
set.
\end{corollary}
\noindent{\bf Proof.} According to Lemma \ref{lem} and Remark
\ref{att}, without loss of generality, we can assume that a forward
invariant random open set $U(\omega)$ is a random pre-attractor
corresponding to the attractor $A(\omega)$. Hence by Lemma 4.2 of
\cite{Cra} the corollary follows. \hfill$\square$

\begin{remark}\label{rem}\rm
(i) By Lemma \ref{lem} and Remark \ref{att} we know that, for
arbitrary random attractor $A(\omega)$, there exists a forward
invariant random pre-attractor $U(\omega)$ which determines
$A(\omega)$. Moreover, by the fact $\Omega_U(\omega)=A(\omega)$, it
follows that $A(\omega)$ pull-back attracts $\overline{U(\omega)}$.
Hence it pull-back attracts any random closed set inside
$U(\omega)$. The random attractor defined in \cite{Cra} only weakly
attracts (i.e. attracts in probability) all random closed sets
inside its fundamental neighborhood. Therefore, the random attractor
we defined in \cite{Liu} (hence in present paper) is stronger than
that of \cite{Cra}. On the other hand, the attractor defined in
\cite{Liu} and \cite{Liu1} are equivalent. In fact, it is clear that
the interior of fundamental neighborhood \cite{Liu1} plays the role
of random pre-attractor \cite{Liu}; conversely, the closure of
random
pre-attractor plays the role of fundamental neighborhood.\\
(ii) For a given random local attractor $A(\omega)$, we can define
the {\em random repeller corresponding to $A(\omega)$} to be
$R(\omega):=X-B(A)(\omega)$ similar to that in \cite{Cra,Liu1} and
we call the pair $(A,R)$ a {\em random attractor-repeller pair}. In
fact, $R(\omega)$ is a random repeller with respect to its random
neighborhood $X-\overline{U(\omega)}$, see \cite{Liu1} or
forthcoming Lemma \ref{non} for the proof. As for the properties of
random attractor-repeller pair, the reader can refer to
\cite{Cra,Liu1} for details.
\end{remark}

Therefore, the main theorem in \cite{Liu} can be restated in the
following theorem, which characterizes random chain recurrent set
using its original form in the deterministic case.

\begin{theorem}\label{thm}
Assume $X$ is a compact metric space, $A(\omega)$ is a random local
attractor and $R(\omega)$ is the random repeller corresponding to
$A(\omega)$, then
\begin{equation*}%\label{mr2}
\mathcal{CR}_\varphi(\omega)=\bigcap [A(\omega)\cup R(\omega)]
\end{equation*}
almost surely, where the intersection is taken over all random local
attractors.
\end{theorem}

\section{Complete Lyapunov function}

\begin{lemma}\label{lya1}
Assume $(A,R)$ is an attractor-repeller pair of $\varphi$, then
there exists an $\mathcal F\times\mathcal B(X)$-measurable Lyapunov function $l$ for $(A,R)$ such that:\\
(i) $l(\omega,x)=0$ when $x\in A(\omega)$, and $l(\omega,x)=1$
     when $x\in R(\omega)$;\\
(ii) for $\forall x\in X\backslash(A(\omega)\cup R(\omega))$ and
$\forall t>0$,
$1>l(\omega,x)>l(\theta_t\omega,\varphi(t,\omega)x)>0$.
\end{lemma}
\noindent{\bf Proof.} The idea of the proof is originated from
\cite{BS,Ar2}. Assume $U(\omega)$ is a forward invariant random
pre-attractor of $A(\omega)$, and we define the first entrance time
of $\varphi(t,\omega)x$ into $\overline{U(\theta_t\omega)}$ as
follows:
\begin{equation}\label{tau}
\tau(\omega,x):=\left\{
\begin{array}{ll}
    -\infty, & x\in A(\omega); \\
    \inf\{t\in\mathbb R|~\varphi(t,\omega)x\in \overline{U(\theta_t\omega)}\},
      &  x\in X\backslash(A(\omega)\cup R(\omega));\\
    +\infty, & x\in R(\omega). \\
\end{array}
\right.
\end{equation}

Since $\omega\mapsto d(x,\overline{U(\omega)})$ is measurable,
$x\mapsto d(x,\overline{U(\omega)})$ is continuous, we have that
$(\omega,x)\mapsto d(x,\overline{U(\omega)})$ is measurable. Hence
for arbitrary $t\in\mathbb R$, $(\omega,x)\mapsto
d(\varphi(t,\omega)x,\overline{U(\theta_t\omega)})$ is measurable.
For $\forall a\in\mathbb R$, we have
\[
\{(\omega,x)|~\tau(\omega,x)\ge a\}=\bigcap_{t<a,t\in\mathbb
Q}\{(\omega,x)|~d(\varphi(t,\omega)x,\overline{U(\theta_t\omega)})>0\},
\]
which verifies that $(\omega,x)\mapsto \tau(\omega,x)$ is
measurable.

By the definition of $\tau(\omega,x)$ in (\ref{tau}), we have
\begin{align*}
\tau(\theta_t\omega,\varphi(t,\omega)x)&=\inf\{s\in\mathbb
R|~\varphi(s,\theta_t\omega)\circ\varphi(t,\omega)x\in\overline{U(\theta_{t+s}\omega)}\}\\
&=\inf\{s\in\mathbb
R|~\varphi(t+s,\omega)x\in\overline{U(\theta_{t+s}\omega)}\}\\
&=\tau(\omega,x)-t.
\end{align*}
Define
\[
l(\omega,x)=\left\{
\begin{array}{ll}
    \frac12 {\rm e}^{\tau(\omega,x)}, & -\infty\le\tau(\omega,x)<0; \\
    \frac12(1+\frac2\pi\arctan\tau(\omega,x)), & 0\le\tau(\omega,x)\le +\infty. \\
\end{array}
\right.
\]
Since $\tau(\omega,x)$ is $\mathcal F\times\mathcal
B(X)$-measurable, hence $l(\omega,x)$ is. It is obvious that the so
defined $l(\omega,x)$ satisfies (i), while (ii) follows from the
fact $\tau(\theta_t\omega,\varphi(t,\omega)x)=\tau(\omega,x)-t$.
This terminates the proof of the lemma. \hfill$\Box$

Denote $\tilde{\mathcal M}=\{x|~x: \Omega\mapsto X~ {\rm
measurable}\}$. For any $x(\omega),y(\omega)\in\tilde{\mathcal M}$,
define
\[d_{\tilde{\mathcal M}}(x,y)=\int_\Omega d_X(x(\omega),y(\omega)){\rm d}\omega.\]
It is clear that $d_{\tilde{\mathcal M}}(\cdot,\cdot)$ is a metric
on $\tilde{\mathcal M}$. By the separability of $\Omega$ and $X$, we
can assume that $\{U_i\}_{i=1}^\infty$ is a countable basis for the
topology of $\Omega$, and $\{x_j\}_{j=1}^\infty$ are countable dense
subsets of $X$. For $U_i,i=1,\cdots,n$ with $U_i\cap
U_j=\emptyset,i\neq j$ and $x_k,k=1,\cdots,n+1$, we define the
following step-like function:
\[
\chi_{U_1,\cdots,U_n}(x_1,\cdots,x_{n+1})=\left\{\begin{array}{ll}
 x_1,&~ \omega\in U_1, \\
\cdots &~ \cdots\\
x_n, &~ \omega\in U_n, \\
x_{n+1}, &~ \omega\in\Omega-\bigcup_{i=1}^n U_i.\\
\end{array}
\right.
\]
Then it is easy to see that the step-like functions defined above
are countable and they constitute a countable dense subset of
$\tilde{\mathcal M}$. For definiteness, we denote this countable
dense subset of $\tilde{\mathcal M}$ as $\mathcal
M=\{x_n(\omega)\}_{n\in\mathbb N}$. Similarly, denote $\mathcal N$
the countable dense subset of $\tilde{\mathcal N}$, the set of all
measurable maps from $\Omega$ to $\mathbb R^+$. By the proof of
Lemma 3.7 in \cite{Liu}, for any given $x(\omega)\in\mathcal
M,\epsilon(\omega),T(\omega)\in\mathcal N$, we can obtain a random
pre-attractor as expressed by (3.11) in \cite{Liu}, and we denote it
here by $U(x,\epsilon,T)$. Therefore we can obtain a random local
attractor $A(x,\epsilon,T)$ determined by $U(x,\epsilon,T)$. Note
that the random local attractors obtained in this way are countable,
and we denote all of them by $\mathcal A$.

\begin{lemma}\label{coun}
Assume $(\Omega,\mathcal F,\mathbb P)$ is a separable metric space
endowed with a probability measure, where $\mathcal F=\mathcal
B(\Omega)$, then
\begin{equation}\label{count}
\mathcal{CR}_\varphi(\omega)=\bigcap_{A\in\mathcal A}[A(\omega)\cup
R(\omega)]
\end{equation}
almost surely, where $\mathcal A$ is defined as above.
\end{lemma}
\noindent{\bf Proof.} By Theorem \ref{thm} we only need to prove the
right hand side of (\ref{count}) is contained in the left hand side.
Note that $\mathcal{N,M}$ is dense in
$\tilde{\mathcal{N}},\tilde{\mathcal{M}}$, respectively. The
remaining proof is similar to that of Lemma 3.7 in \cite{Liu}, here
we omit the details.  \hfill $\Box$

\begin{definition}\label{LY}
A {\em complete Lyapunov function} for the RDS $\varphi$ is an
$\mathcal F\times\mathcal B(X)$-measurable function $L:\Omega\times
X\mapsto\mathbb R^+$ with the following properties:\\
(1) $L(\theta_t\omega,\varphi(t,\omega)x)=L(\omega,x)$,
               $\forall t>0$ when $x\in\mathcal{CR}_\varphi(\omega)$;\\
(2) If $x\in X-\mathcal{CR}_\varphi(\omega)$, then
          \[
          L(\theta_t\omega,\varphi(t,\omega)x)<L(\omega,x),~\forall
          t>0;
          \]
(3) The range of $L$ on $\mathcal{CR}_\varphi(\omega)$ is a compact
nowhere dense
               subset of $[0,1]$;\\
(4) $L$ separates different random chain transitive components of
$\varphi$ (to be defined below).
\end{definition}

For definiteness, denote the above $\mathcal A=\{A_n(\omega)|~n\in
\mathbb N\}$ and assume $l_n(\omega,x)$ is the Lyapunov function
determined by the random attractor-repeller $(A_n,R_n)$ in Lemma
\ref{lya1}. Define
\begin{equation}\label{L}
L(\omega,x)=\sum_{n=1}^{\infty}\frac{2l_n(\omega,x)}{3^n}.
\end{equation}
Later we will show that $L(\omega,x)$ defined by (\ref{L}) is a
complete Lyapunov function for $\varphi$. Firstly it is easy to see
that $L(\omega,x)$ is well defined since the sum in (\ref{L}) is
uniformly convergent and $L(\omega,x)$ is $\mathcal F\times\mathcal
B(X)$-measurable since each $l_n(\omega,x)$ is.

A {\em critical value of $L$} is defined to be one achieved on the
random chain recurrent set. For any $A(\omega)\in\mathcal A$, we
have $\mathcal{CR}_\varphi(\omega)\subset A(\omega)\cup R(\omega)$
by Lemma \ref{coun}. Hence each $l_n(\omega,x)$ is either $0$ or $1$
at a point of the random chain recurrent set. Therefore the set of
the critical values of $L$ is a compact nowhere dense subset of
$[0,1]$. For a given critical value $c$ of $L$, denote
$L^{-1}(c)=\Omega\times\{\tilde C(\omega)\}$. Then by the
measurability of $L$, the graph of $\tilde C(\omega)$
\[
{\rm graph}(\tilde C)=\{(\omega,x)|~x\in\tilde
C(\omega)\}=\{(\omega,x)|~L(\omega,x)=c\}
\]
is $\mathcal F\times\mathcal B(X)$-measurable. Hence $\tilde
C(\omega)$ is an $\mathcal F^u$-measurable random closed set by
Theorem III.30 on page 80 of \cite{Cas}. Denote $C(\omega)=\tilde
C(\omega)\cap\mathcal{CR}_\varphi(\omega)$, then $C(\omega)$ is an
$\mathcal F^u$-measurable random compact set (since
$\mathcal{CR}_\varphi(\omega)$ is an $\mathcal F^u$-measurable
random compact set). For convenience, we also use
$L^{-1}(c)\cap\mathcal{CR}_\varphi(\omega)$ to denote $C(\omega)$,
indicating that $C(\omega)$ is determined by the critical value $c$
of $L$. Since $C(\omega)$ is determined by the critical value $c$,
we have
\begin{equation}\label{ddis}
C(\omega)\subset A(\omega) ~\mathbb P{\rm -a.s.} ~{\rm or}~
C(\omega)\subset R(\omega) ~\mathbb P{\rm -a.s.}
\end{equation}
for any $A(\omega)\in\mathcal A$, where $(A,R)$ is a random
attractor-repeller pair.

We call $C(\omega)$ obtained in this way a {\em random chain
transitive component of $\varphi$}. By Lemma \ref{coun} we have that
a chain transitive component $C(\omega)$ can be expressed by
\begin{equation}\label{component}
C(\omega)=\bigcap_{n\in\mathbb N}C_n(\omega),
\end{equation}
where $C_n(\omega)=A_n(\omega)$ or $C_n(\omega)=R_n(\omega)$,
recalling that we have assumed that $\mathcal A=\{A_n(\omega)|~n\in
\mathbb N\}$.

\begin{definition}
Assume $C,C'$ are distinct random chain transitive components of
$\varphi$. We say that a random local attractor $A(\omega)$ {\em
distinguishes} between them if either $C(\omega)\subset A(\omega)$
and $C'(\omega)\subset R(\omega)$ or if the corresponding conditions
hold with $C$ and $C'$ interchanged.
\end{definition}

An equivalence relation on the set of random chain recurrent
variables is given
by:\\
\quad $x(\omega)\sim y(\omega)$ if and only if for arbitrary
$\epsilon(\omega),T(\omega)>0$, there is one
$\epsilon(\omega)$-$T(\omega)$-chain from $x(\omega)$ to $y(\omega)$
and another one from $y(\omega)$ to $x(\omega)$ $\mathbb P$-a.s.

With respect to the relation between random chain transitive
component of $\varphi$ and the equivalence relation introduced
above, we have the following proposition:

\begin{proposition}\label{compo}
(1) Assume $C(\omega)$ is a random chain transitive component of
$\varphi$, then for $\forall x(\omega),y(\omega)\in C(\omega)$, we
have $x(\omega)\sim y(\omega)$. \\
(2) Assume $C'(\omega)$ is another random chain transitive component
of $\varphi$, then for $\forall x(\omega)\in C(\omega)$ and $\forall
z(\omega)\in C'(\omega)$, we have $x(\omega)\not\sim z(\omega)$.
\end{proposition}
\noindent{\bf Proof.} (1) We only need to prove that for
$\forall\epsilon(\omega),T(\omega)\in\mathcal N$, there exist
$\epsilon(\omega)$-$T(\omega)$-chains from $x(\omega)$ to
$y(\omega)$ and from $y(\omega)$ to $x(\omega)$ $\mathbb P$-a.s. If
this is false, say, there exist $\epsilon_0(\omega),
T_0(\omega)\in\mathcal N$ such that there is no
$2\epsilon_0$-$T_0$-chain from $x(\omega)$ to $y(\omega)$ $\mathbb
P$-a.s. By the density of $\mathcal M$ in $\tilde{\mathcal M}$,
there exists some $x_1(\omega)\in\mathcal M$ such that $x(\omega)\in
B(A(x_1,\epsilon_0,T_0))$ $\mathbb P$-a.s. Since $x(\omega)$ is
random chain recurrent, we have $x(\omega)\in A(x_1,\epsilon_0,T_0)$
$\mathbb P$-a.s. by Lemma 3.6 of \cite{Liu}. Since there is no
$2\epsilon_0$-$T_0$-chain from $x(\omega)$ to $y(\omega)$ $\mathbb
P$-a.s., we can obtain that $y(\omega)\notin A(x_1,\epsilon_0,T_0)$
with positive probability. In fact, if $y(\omega)\in
A(x_1,\epsilon_0,T_0)$ $\mathbb P$-a.s., then we have $y(\omega)\in
A(x_1,\epsilon_0,T_0)\subset U(x,2\epsilon_0,T_0)$ $\mathbb P$-a.s.
when $x_1(\omega)$ is close enough to $x(\omega)$. Hence we obtain
that there is $2\epsilon_0$-$T_0$-chain from $x(\omega)$ to
$y(\omega)$ $\mathbb P$-a.s. by the definition of
$U(x,2\epsilon_0,T_0)$, a contradiction. \\
(2) To see this, notice that
\[
\mathcal{CR}_\varphi(\omega)=\bigcup_{c~{\rm
critical~value}}[L^{-1}(c)\cap\mathcal{CR}_\varphi(\omega)].
\]
Assume that the two random transitive components $C(\omega)$ and
$C'(\omega)$ are determined by two critical values $c\neq c'$ of
$L$, respectively. That is
\[
C(\omega)=L^{-1}(c)\cap\mathcal{CR}_\varphi(\omega),~~
C'(\omega)=L^{-1}(c')\cap\mathcal{CR}_\varphi(\omega).
\]
Since $c\neq c'$, there exists some $n_0\in\mathbb N$ such that
$l_{n_0}$ is 0 on $C(\omega)$ and is 1 on $C'(\omega)$ or the
converse holds. Hence we have
\[
C(\omega)\subset A_{n_0}(\omega),\quad C'(\omega)\subset
R_{n_0}(\omega)~~\mathbb P{\rm -a.s.}
\]
or
\[
C'(\omega)\subset A_{n_0}(\omega),\quad C(\omega)\subset
R_{n_0}(\omega)~~\mathbb P{\rm -a.s.}
\]
Without loss of generality, we assume the former case holds. If for
any $\epsilon(\omega),T(\omega)>0$, there exists an
$\epsilon(\omega)$-$T(\omega)$-chain from $x(\omega)$ to $z(\omega)$
$\mathbb P$-a.s. Then by the proof of Lemma 3.4 of \cite{Liu}, we
obtain that $z(\omega)\in A_{n_0}(\omega)$ $\mathbb P$-a.s., a
contradiction to $z(\omega)\in C'(\omega)$. That is
$x(\omega)\not\sim z(\omega)$. \hfill $\Box$

\begin{remark}\label{comre}\rm
By the above proof we know that if $C,C'$ are distinct random chain
transitive components of $\varphi$, then there exists an element
$A(\omega)\in\mathcal A$ that distinguishes between them.
\end{remark}

\begin{lemma}\label{inclu}
Assume $x(\omega),y(\omega)\in\mathcal{CR}_\varphi(\omega)$ $\mathbb
P$-a.s. with the property that for $\forall\epsilon(\omega)$,
$T(\omega)>0$ there is a random $\epsilon(\omega)$-$T(\omega)$-chain
from $x(\omega)$ to $y(\omega)$ $\mathbb P$-a.s. If $A(\omega)$ is a
random local attractor containing $x(\omega)$ $\mathbb P$-a.s., then
$A(\omega)$ also contains $y(\omega)$ $\mathbb P$-a.s.
\end{lemma}
\noindent{\bf Proof.} If $y(\omega)\notin A(\omega)$ with positive
probability, then by Lemma 3.6 of \cite{Liu} we have $y(\omega)\in
X-B(A)(\omega)$ with positive probability. Since for
$\forall\epsilon(\omega),T(\omega)>0$, there is a random
$\epsilon(\omega)$-$T(\omega)$-chain from $x(\omega)$ to $y(\omega)$
$\mathbb P$-a.s., similar to the proof of Lemma 3.4 of \cite{Liu}
and by the fact $x(\omega)\in A(\omega)$, we obtain that
$y(\omega)\in U(\omega)\subset B(A)(\omega)$ $\mathbb P$-a.s., where
$U(\omega)$ is a random pre-attractor which determines $A(\omega)$,
a contradiction. \hfill $\Box$

\begin{corollary}\label{col2}
Assume $C,C'$ are two random chain transitive components of
$\varphi$ with the property that for $\forall\epsilon(\omega)$,
$T(\omega)>0$ there is a random $\epsilon(\omega)$-$T(\omega)$-chain
from $C$ to $C'$ $\mathbb P$-a.s. If $A(\omega)$ is a random local
attractor containing $C$ $\mathbb P$-a.s., then $A(\omega)$ also
contains $C'$ $\mathbb P$-a.s.
\end{corollary}
\noindent{\bf Proof.} The corollary follows directly from Lemma
\ref{inclu}. \hfill $\Box$

Now we are ready to prove our main result--Theorem \ref{main}.

\begin{theorem}\label{smain} {\bf (complete Lyapunov function).}
The function defined by (\ref{L}) is a complete Lyapunov function
for the RDS $\varphi$. Moreover, if $C$ and $C'$ are distinct random
chain transitive components of $\varphi$ with the property that for
arbitrary $\epsilon(\omega),T(\omega)>0$ there is an
$\epsilon(\omega)$-$T(\omega)$-chain from $C$ to $C'$ $\mathbb
P$-a.s., then $L(\Omega,C)>L(\Omega,C')$.
\end{theorem}
\noindent{\bf Proof.} We first verify (1)--(4) of Definition \ref{LY} one by one.\\
(1): This follows from the fact that if
$x\in\mathcal{CR}_\varphi(\omega)$, then
\[
l_n(\omega,x)= l_n(\theta_t\omega,\varphi(t,\omega)x),\forall t>0
\]
takes value $0$ or $1$ for each $n\in\mathbb N$.\\
(2): Since $x\in X-\mathcal{CR}_\varphi(\omega)$, there exists an
$A_n\in\mathcal A$ such that $x\in X-A_n(\omega)\cup R_n(\omega)$.
Then by Lemma \ref{lya1} we have
\[
l_n(\theta_t\omega,\varphi(t,\omega)x)< l_n(\omega,x), ~\forall t>0.
\]
(3): It has been proved above.\\
(4): By the definition of random chain transitive component, it is
clear that $L$ is constant on each random chain transitive component
and $L$ takes different values on different random chain transitive
components. This verifies (4).

By Corollary \ref{col2}, any random local attractor containing $C$
must also contain $C'$, therefore $l_n=0$ on $C$ implies $l_n=0$ on
$C'$, that is $l_n(\Omega,C)\ge l_n(\Omega,C')$ for each
$n\in\mathbb N$. Hence we have obtained $L(\Omega,C)\ge
L(\Omega,C')$, this together with (4) verifies
$L(\Omega,C)>L(\Omega,C')$.\hfill$\square$

It should be pointed out that random chain transitive components are
not determined by the complete Lyapunov function, but by the RDS
$\varphi$ itself. See (\ref{component}).

Now we give two simple examples to illustrate our results.
\begin{example}\rm
(1) In Example 4.1 of \cite{Liu}, we have shown that
$\mathcal{CR}_\varphi(\omega)=\{-1,0,1\}$ $\mathbb P$-a.s. By
Theorem \ref{main} we can construct a complete Lyapunov function for
$\varphi$. It is clear that $C_1=\{-1\}$, $C_2=\{0\}$, $C_3=\{1\}$
are three random chain transitive components of $\varphi$. Here
$W_t$ is a Wiener process for which we take two-sided time
$t\in\mathbb R$. This is standard to put a stochastic differential
equation in the framework of RDS, see \cite{Ar1,Ar2,Chu,Cra1}
etc for details. \\
(2) In Example 4.2 of \cite{Liu}, we have shown that the RDS
$\varphi$ has no non-trivial attractor besides $X$ and $\emptyset$.
Hence $\mathcal{CR}_\varphi(\omega)=X$ $\mathbb P$-a.s. Therefore,
by (\ref{component}) we obtain that the only random chain transitive
component is $X$, i.e. the random chain recurrent set itself.
\hfill$\Box$
\end{example}

\section{Extension to noncompact spaces}

We know that assuming $X$ being compact is too restrictive for
applications, so in this section we assume that $X$ is a Polish
space and try to extend Conley decomposition theorem for RDS to
noncompact Polish spaces.

\begin{definition}\label{absor}
Assume $\epsilon(\omega)>0$ is a random variable. \\
(i) A random open set $U(\omega)$ is called {\em
$\epsilon$-absorbing} if there exists a random variable
$T(\omega)>0$ such that $U(\omega)$ contains the
$\epsilon$-neighborhood of $U_T(\omega):=\overline{\bigcup_{t\ge
T}\varphi(t,\theta_{-t}\omega)U(\theta_{-t}\omega)}$, i.e.
\[
B_{\epsilon}(U_T(\omega))\subset U(\omega).
\]
And we call a random open set $U(\omega)$ {\em absorbing} if it is
$\epsilon$-absorbing for some random variable
$\epsilon(\omega)>0$.\\
(ii) A random open set $V(\omega)$ is called {\em
$\epsilon$-repelling} if there exists a random variable
$T(\omega)>0$ such that $V(\omega)$ contains the
$\epsilon$-neighborhood of $\hat
V_T(\omega):=\overline{\bigcup_{t\ge
T}\varphi(-t,\theta_{t}\omega)V(\theta_{t}\omega)}$, i.e.
\[
B_{\epsilon}(\hat V_T(\omega))\subset V(\omega).
\]
And we call a random open set $V(\omega)$ {\em repelling} if it is
$\epsilon$-repelling for some random variable $\epsilon(\omega)>0$.
\end{definition}

\begin{definition}\label{def}
(i) An invariant random closed set $A(\omega)$ is called an {\em
(local) attractor} if there exists an absorbing neighborhood
$U(\omega)$ of $A(\omega)$ such that $A(\omega)=\Omega_U(\omega)$.
And we call
\[
B(A,U)(\omega):=\{x|~\varphi(t,\omega)x\in U(\theta_t\omega)~{\rm
for~ some}~t\ge 0\}
\]
the {\em basin of attraction of $A(\omega)$ with respect to
$U(\omega)$}. \\
(ii) An invariant random closed set $R(\omega)$ is called a {\em
(local) repeller} if there exists a repelling neighborhood
$V(\omega)$ of $R(\omega)$ such that $R(\omega)=\alpha_U(\omega)$.
And we call
\[
B(R,V)(\omega):=\{x|~\varphi(t,\omega)x\in V(\theta_t\omega)~{\rm
for~ some}~t\le 0\}
\]
the {\em basin of repulsion of $R(\omega)$ with respect to
$V(\omega)$}.
\end{definition}

\begin{remark}\rm (i) When $X$ is compact, the basin of attraction of an
attractor is independent of the choice of absorbing neighborhoods;
but when $X$ is not compact, the basin is generally dependent on
absorbing neighborhoods. Of course the same conclusion holds for
basin of repulsion.
See \cite{Hu0} for details.\\
(ii) For any random set $D(\omega)$, $\Omega_D(\omega)$ is invariant
even if $X$ loses compactness. See Lemma 3.2 and Remark 3.7 of
\cite{Cra1} for details. Similarly, the alpha-limit set
$\alpha_D(\omega)$ is also invariant.
\end{remark}

\begin{lemma}\label{non}
Assume $\varphi$ is an RDS on a Polish space $X$, and assume
$A(\omega)$ is an attractor of $\varphi$ with a forward invariant
absorbing neighborhood $U(\omega)$ and the basin of attraction with
respect to $U$, $B(A,U)(\omega)$. Then $R(\omega):=X-B(A,U)(\omega)$
is a random repeller with a repelling neighborhood
$V(\omega):=X-\overline{U(\omega)}$ and the basin of repulsion with
respect to $V(\omega)$ being $B(R,V)(\omega)=X-A(\omega)$. (We call
this $R(\omega)$ {\em the repeller corresponding to $A(\omega)$ with
respect to $U(\omega)$}.)
\end{lemma}
\noindent{\bf Proof.} The proof is similar to that of Lemma 4.1 in
\cite{Liu1}. By the forward invariance of $U(\omega)$, we have
\[
\varphi(t,\theta_{-t}\omega)U(\theta_{-t}\omega)\subset
\varphi(s,\theta_{-s}\omega)U(\theta_{-s}\omega),~{\rm for}~\forall
t\ge s, s,t\in\mathbb R.
\]
Hence we have
\begin{align}\label{5.1}
A(\omega)&=\bigcap_{T\ge 0}\overline{\bigcup_{t\ge
T}\phi(t,\theta_{-t}\omega)U(\theta_{-t}\omega)}=\bigcap_{n\in\mathbb
N}\varphi(n,\theta_{-n}\omega)\overline{U(\theta_{-n}\omega)}\nonumber\\
&=\lim_{n\rightarrow\infty}\varphi(n,\theta_{-n}\omega)\overline{U(\theta_{-n}\omega)}
\end{align}
and
\begin{align}\label{basin}
B(A,U)(\omega)&=\bigcup_{t\ge
0}\varphi(-t,\theta_t\omega)U(\theta_t\omega)=\bigcup_{n\in\mathbb
N}\varphi(-n,\theta_n\omega)U(\theta_n\omega)\nonumber\\
&=\lim_{n\rightarrow\infty}\varphi(-n,\theta_{n}\omega)U(\theta_{n}\omega).
\end{align}

Since $U(\omega)$ is a forward invariant random open set, we have
that $V(\omega):=X-{\overline{U(\omega)}}$ is a backward invariant
random open set (see page 35 of \cite{Ar1}). Denote $\tilde
R(\omega):=\alpha_{V}(\omega)$, then $\tilde R(\omega)$ is a random
repeller with a repelling neighborhood $V(\omega)$. By the
definition of alpha-limit set, we have $R(\omega)\subset\tilde
R(\omega)$ due to the fact $R(\omega)\subset V(\omega)$ and the
invariance of $R(\omega)$. If there exists some $x_0\in\tilde
R(\omega)\backslash R(\omega)$, then $x_0\in B(A,U)(\omega)$.
Therefore there exists some $t_0\ge 0$ such that
$\varphi(t_0,\omega)x_0\in U(\theta_{t_0}\omega)$. Noticing that
$\tilde R(\omega)$ is an invariant random closed set, we have
$\varphi(t_0,\omega)x_0\in\tilde R(\theta_{t_0}\omega)$. This is a
contradiction to the fact $\tilde R(\omega)\cap U(\omega)=\emptyset$
for each $\omega$. Therefore we obtain $R(\omega)=\tilde R(\omega)$,
i.e. $R(\omega)$ is a repeller with a repelling neighborhood
$V(\omega)$. Now we will show $B(R,V)(\omega)=X-A(\omega)$. In fact
\begin{align}
B(R,V)(\omega)&=\bigcup_{n\in\mathbb
N}\varphi(n,\theta_{-n}\omega)V(\theta_{-n}\omega)\label{1}\\
&=\lim_{n\rightarrow\infty}\varphi(n,\theta_{-n}\omega)V(\theta_{-n}\omega)\label{2}\\
&=\lim_{n\rightarrow\infty}\varphi(n,\theta_{-n}\omega)[X-\overline{U(\theta_{-n}\omega)}]\nonumber\\
&=\lim_{n\rightarrow\infty}[X-\varphi(n,\theta_{-n}\omega)\overline{U(\theta_{-n}\omega)}]\label{3}\\
&=X-A(\omega),\label{4}
\end{align}
where (\ref{1}) and (\ref{2}) hold completely similar to
(\ref{basin}) if we take $t=-t$, (\ref{3}) follows from the fact
that $\varphi(n,\omega)$ is a homeomorphism on $X$ and (\ref{4})
holds by (\ref{5.1}). This completes the proof of the
lemma.\hfill$\square$

Similar to the proofs of Lemmas 3.4, 3.6 and 3.7 in \cite{Liu}, we
can obtain the following result (Notice that in these proofs, the
compactness is only relevant in proving Lemma 3.4, where we use
compactness to conclude $\bar d(\omega)>0$. In Definitions
\ref{absor} and \ref{def}, by requesting that the neighborhood of
attractor be absorbing we overcome this difficulty when $X$ loses
compactness.):

\begin{theorem}\label{th}
Assume $X$ is a Polish space, $U(\omega)$ is an absorbing set,
$A(\omega)$ is the random local attractor determined by $U(\omega)$,
and $B(A,U)(\omega)$ is the basin of $A(\omega)$ with respect to
$U(\omega)$, then
\begin{equation*}%\label{mr}
X-\mathcal{CR}_\varphi(\omega)=\bigcup [B(A,U)(\omega)-A(\omega)]
\end{equation*}
almost surely, where the union is taken over all absorbing sets.
\end{theorem}

\begin{lemma}\label{lem1}
Assume $A(\omega)$ is an attractor with an absorbing neighborhood
$U(\omega)$ which determines $A(\omega)$. Then there is a forward
invariant absorbing neighborhood $\tilde U(\omega)$ of $A(\omega)$
such that $\Omega_{\tilde U}(\omega)=A(\omega)$.
\end{lemma}
\noindent{\bf Proof.} Noting that the compactness is not required
when we prove $\Omega_{\tilde
U}(\omega)=\Omega_{U}(\omega)=A(\omega)$ in Lemma \ref{lem}. Hence
there exists a forward invariant neighborhood $\tilde U(\omega)$ of
$A(\omega)$ such that $\Omega_{\tilde U}(\omega)=A(\omega)$ when $X$
loses compactness. We only need to show that $\tilde U(\omega)$ is
absorbing. Assume that there exist $\epsilon(\omega), T(\omega)>0$
such that
\[
B_{\epsilon}(U_T(\omega))\subset U(\omega).
\]
By the fact $\Omega_{\tilde U}(\omega)=\Omega_{U}(\omega)=A(\omega)$
we obtain that there exists $t=t(\omega)$ such that
\[
\tilde U_t(\omega)\subset U_T(\omega),
\]
where $\tilde U_t(\omega)$ is defined similar to $U_T(\omega)$.
Define
\[
\tilde T(\omega):=\inf\{t\in\mathbb
R^+|~\varphi(t,\theta_{-t}\omega)\overline{\tilde
U(\theta_{-t}\omega)}\subset U_T(\omega)\},
\]
then $\tilde T(\omega)$ is measurable by similar argument to the
proof of Lemma 3.5 in \cite{Liu}. Thus $B_{\epsilon}(\tilde
U_{\tilde T}(\omega))\subset U(\omega)\subset\tilde U(\omega)$, and
hence $\tilde U(\omega)$ is absorbing. The proof is complete.
\hfill$\square$

Similar to the case when $X$ is compact, we have the following
theorem:

\begin{theorem}\label{thm1}
Assume $X$ is a Polish space, $U(\omega)$ is a forward invariant
absorbing set, $A(\omega)$ is the random local attractor determined
by $U(\omega)$ and $R(\omega)$ is the repeller corresponding to
$A(\omega)$ with respect to $U(\omega)$, then
\begin{equation*}%\label{mr1}
\mathcal{CR}_\varphi(\omega)=\bigcap [A(\omega)\cup R(\omega)]
\end{equation*}
almost surely, where the intersection is taken over all {\em forward
invariant} absorbing sets.
\end{theorem}
\noindent{\bf Proof.} By Theorem \ref{th}, we only need to verify
\begin{equation}\label{eq}
\bigcap_{U\in\mathcal U} [A(\omega)\cup
R(\omega)]\subset\mathcal{CR}_\varphi(\omega),
\end{equation}
where $\mathcal U$ denotes the set of all forward invariant
absorbing sets. To this end, we only need to show
\begin{equation}\label{eq1}
\bigcup [B(A,U)(\omega)-A(\omega)]\subset \bigcup_{U\in\mathcal U}
[B(A,U)(\omega)-A(\omega)],
\end{equation}
where the union on the left hand is taken over all absorbing sets.
By Lemma \ref{lem1} we know that for an attractor $A(\omega)$ with
an absorbing neighborhood $U(\omega)$, there is an absorbing
neighborhood $\tilde U(\omega)\in\mathcal U$ of $A(\omega)$ such
that $\Omega_{\tilde U}(\omega)=\Omega_{U}(\omega)=A(\omega)$ and
$U(\omega)\subset\tilde U(\omega)$. Thus we have
\[
B(A,U)(\omega)\subset B(A,\tilde U)(\omega).
\]
Therefore we have proved (\ref{eq1}) and hence the theorem.
\hfill$\square$

By mimicking the proof when the state space is compact, we can
obtain the complete Lyapunov function for $\varphi$ when $X$ is not
compact; furthermore, we can discuss the chain transitive components
etc completely similar to the compact case. In fact, the compactness
is not relevant during these steps. Hence we omit details here.
Therefore we have the random Conley decomposition theorem on Polish
spaces:

\begin{theorem}{\bf (random Conley decomposition theorem on Polish spaces).}
Any random dynamical system (with a separable metric space endowed
with a probability measure as base space) on a Polish space
decomposes the Polish space into a random chain recurrent part and a
random gradient-like part.
\end{theorem}

\begin{example}\rm
Consider the Lorenz system in $\mathbb R^3$ described by the
equations:
\[
\left\{
  \begin{array}{l}
    \dot x=\sigma(y-x),\\
    \dot y=\rho x-y-xz, \\
    \dot z=xy-\beta z
  \end{array}
\right.
\]
with parameters $\sigma,\rho,\beta>0$. Assume that the parameter
$\sigma,\rho,\beta$ are perturbed by real noises, i.e.
\begin{eqnarray*}
&\sigma(\omega)=\sigma+\xi(\omega),\\
&\rho(\omega)=\rho+\eta(\omega),\\
&\beta(\omega)=\beta+\zeta(\omega).
\end{eqnarray*}
Assume that the perturbed parameters are still positive and
$\rho(\omega)<\sigma(\omega)\le 1$ almost surely. Then the perturbed
random ODE generates an RDS and we denote it by $\varphi$, see
\cite{Ar1} for details. Let
\[
L(\omega,x,y,z)=x^2+y^2+z^2
\]
and define the random set
\[
D_r(\omega):=\{(x,y,z)|~L(\omega,x,y,z)\le r\}.
\]
Then clearly $D_r(\omega)$ is a random compact set. Moreover, the
derivative of $L$ with respect to $t$ along the orbits of $\varphi$
\begin{align}\label{exam}
\frac{{\rm d}L}{{\rm d}t}&=2x(\sigma(y-x))+2y(\rho
x-y-xz)+2z(xy-\beta z)\nonumber\\
&=-2\sigma x^2-2y^2-2\beta z^2+2(\rho+\sigma)xy\nonumber\\
&\le -(\sigma-\rho)x^2-(2-\rho-\sigma)y^2-2\beta z^2\nonumber\\
&<0
\end{align}
whenever $(x,y,z)\neq (0,0,0)$ by the assumption
$\rho(\omega)<\sigma(\omega)\le 1$, where
$\rho=\rho(\theta_t\omega)$, $\sigma=\sigma(\theta_t\omega)$,
$\beta=\beta(\theta_t\omega)$ and $X=X(t,\omega)$ with $X:=(x,y,z)$.
Therefore, for $\forall X\in D_r(\omega)$ and $\forall t>0$ we have
$L(\theta_t\omega,\varphi(t,\omega)X)<L(\omega,X)\le r$, i.e.
$\varphi(t,\omega)X\in D_r(\theta_t\omega)$. Hence $D_r(\omega)$ is
a forward invariant random compact set and clearly it determines a
random attractor $A(\omega)=\{0\}$ with the basin
$B(A)(\omega)=\mathbb R^3$ almost surely. And by (\ref{exam}) we
easily know that $A(\omega)$ is the only non-trivial attractor of
$\varphi$. Therefore, by Theorem \ref{thm1} we obtain that the
random chain recurrent set is $\{0\}$ and hence it is the only
random chain transitive component. \hfill$\Box$
\end{example}

{\em Notes to \cite{Liu}:} The definition of random open set in
\cite{Liu} should be same as Definition \ref{set} of present paper.
Correspondingly, the items (i), (ii) and (vi) of Proposition 2.1 in
\cite{Liu} should be stated as follows:
\begin{description}
  \item[(i)] $D(\omega)$ is a random closed set in $X$ if and only if the set $\{\omega :
D(\omega) \bigcap U \neq \emptyset \}$ is measurable for any open
set $U \subset X$;
  \item[(ii)] $D(\omega)$ is measurable if and only
if $\overline{D(\omega)}$ is a random closed set;
  \item[(vi)] If $\{D_n, n \in\mathbb N\}$ is a sequence of random open
  sets and there exists $n_0\in\mathbb N$ such that $D_{n_0}^c$ is a
  random compact set, then $D=\bigcup_{n\in\mathbb N}D_n$ is also a random open
  set. Or if $\{D_n, n \in\mathbb N\}$ is a sequence of measurable
  multifunctions, then $D=\bigcup_{n\in\mathbb N}D_n$ is also a measurable
  multifunction.
\end{description}

By Lemma \ref{lem} and Remark \ref{rem} of present paper, the proofs
of Lemmas 3.1, 3.3 and 3.6 in \cite{Liu} can be greatly simplified
and some minor errors can be easily avoided.

\vskip8mm
\noindent{\bf Acknowledgements} \\
The author expresses his sincere thanks to Professor Yong Li for his
instructions and many invaluable suggestions. The author is very
grateful to Professors Youqing Ji and Xiaoyun Yang as well as
Shuguan Ji and Menglong Su for helpful discussions. Many thanks to
the anonymous referees for their careful reading the manuscript and
invaluable comments which greatly improved the presentation of the
paper.

\end{document}